\documentclass[reqno,12pt,article]{amsart}
\usepackage[cp1251]{inputenc}
\usepackage[russian]{babel}
\usepackage{amsmath,amsthm,amssymb}
\usepackage{multirow}
\begin{document}

\author{Л.\,М.~Арутюнян}
\address{Lomonosov MSU, the faculty of mechanics and mathematics, Moscow city.}
\email{Lavrentin@ya.ru}

\title{Интегральные неравенства Ремеза для многочленов на выпуклых телах}

%\markboth{L. M. Arutyunyan}{On growth of the number of determinants}

\maketitle

\large

\footnotetext{\ Работа поддержана грантом РНФ 17-11-01058 (выполняемым при МГУ
им. М.В.Ломоносова). }

\section{Введение.}

Пусть на некотором пространстве $\Omega$ с вероятностной мерой $\mu$ задана система функций $X$.
Будем называть классическими неравенствами типа Ремеза неравенства вида
$$
\|f\|_{L^{\infty}(\Omega)} \le C(\Omega, \mu(A), X) \|f\|_{L^{\infty}(A)},
$$
которые выполняются для всякого измеримого множества $A$ и всякой функции $f$ из пространства $X$ с константой
$C(\Omega, \mu(A), X)$, которая зависит лишь от пространства $\Omega$, меры множества $A$, а также от класса функций $X$.
Таким образом, классические неравенства типа Ремеза оценивают супремум функции на всем пространстве через ее супремум на некотором
произвольном множестве положительной меры.

Оригинальное неравенство, доказанное Ремезом, относилось к случаю класса $P^d$ многочленов степени $d$
на отрезке $[a,b]$ с мерой $\lambda_{[a,b]}$, которая является вероятностным сужением меры Лебега на отрезок $[a,b]$,
причем в этом случае
$$C(\mathbb{R}, \lambda_{[a,b]}(A), P^d) \le \Bigl(\frac{4}{\lambda_{[a,b]}(A)}\Bigr)^d.$$
Брудный и Гансбург рассмотрели аналогичную задачу в пространстве $\Omega=\mathbb{R}^n$ с многочленами степени $d$, но уже от $n$ переменных, причем
в качестве $\mu$ рассматривалось равномерное распределение $\lambda_V$ на произвольном выпуклом теле $V$ (см. \cite{BrG}). В этом случае константа в неравенстве типа Ремеза
будет стремиться к бесконечности с ростом размерности $n$:
$$C(\mathbb{R}^n, \lambda_{V}(A), P^d) \le \Bigl(\frac{4 n}{\lambda_{[a,b]}(A)}\Bigr)^d.$$

Интегральными неравенствами Ремеза будем называть аналогичные неравенства, но уже для интегральных норм:
$$
\|f\|_{L^{1}(\mu)} \le C(\Omega,\mu(A), X) \|f\|_{L^{1}(\mu_{A})},
$$
где $\mu_A$ --- это нормированное сужение меры $\mu$ на множество $A$, то есть мера $\mu(\cdot) = \frac{\mu(\cdot\, \cap A)}{\mu(A)}$.
Оказывается,
что константа в интегральном неравенстве Ремеза для многочленов
на выпуклых телах может быть взята независимой от размерности,
в частности, такое неравенство можно сформулировать и в бесконечномерном случае, заменив равномерные распределения
на выпуклом компакте на логарифмически вогнутые меры (которые в определенном смысле являются
бесконечномерными аналогами равномерных распределений на выпуклых телах, см. секцию 5.).

\section{Определения и основные результаты.}

Будем рассматривать локально выпуклые пространства $X$ с радоновскими логарифмически вогнутыми мерами на них
(хотя для простоты читатель может считать пространством $\mathbb{R}^n$, а мерой --- равномерное распределение на выпуклом компакте).

Логарифмически вогнутой мерой называется вероятностная мера, для которой при всех борелевских множествах
$A$ и $B$ выполняется неравенство
$$
    \mu(\alpha A + (1-\alpha)B) \ge \mu(A)^{\alpha}\mu(B)^{1-\alpha}\quad \text{ для всяких } \alpha\in [0,1].
$$
Примерами логарифмически вогнутых мер служат равномерные распределения на выпуклых телах и гауссовские меры.

Многочленом степени $d$ на линейном пространстве
$X$ называется функция $f$ вида
$${f(x)=b_0+b_1(x)+\cdots +b_d(x,\ldots,x)},$$
где $b_0$ --- постоянная, $b_k\colon\, X^k\to \mathbb{R}$ --- полилинейная
функция (т.е. линейная по каждому аргументу). Аналогично определяются
полиномиальные степени $d$ отображения из $X$ в линейное пространство~$Y$.

Тригонометрическим полиномом степени $d$ на линейном пространстве $X$ называют
функцию вида
$$f(x)=\sum_{k=1}^d e^{i \ell_k(x)},$$
где $\ell_k$ --- некоторые линейные функционалы.

Если $f$ --- измеримое полиномиальное отображение из пространства $X$ с логарифмически вогнутой мерой $\mu$
в банахово пространство $Y$ с нормой $\|\cdot\|_Y$, то существуют следующие величины (см. \cite{AK, NVS}):
$$
\|f\|_p=\Bigl(\int \|f\|_Y^p d\mu\Bigr)^{1/p},\ \ p\in(-{1}/{d},0)\cup(0,\infty),
$$
$$
\|f\|_0=\lim_{p\rightarrow 0} \|f\|_p.
$$
Кроме того, для краткости введем обозначения
$$
\|f\|_{p,A}=\Bigl(\frac{1}{\mu(A)}\int_{A} \|f\|_Y^p d\mu\Bigr)^{1/p},\ \ p\in(-{1}/{d},0)\cup(0,\infty),
$$
$$
\|f\|_{0,A}=\lim_{p\rightarrow 0} \|f\|_{p,A}.
$$
В дальнейшем все многочлены считаются измеримыми.

Основной результат заключается в следующем.

{\bf Теорема 1.} Пусть $f$ --- полиномиальное отображение степени $d$, мера $\mu$
является логарифмически вогнутой, а множество $A$ имеет положительную меру. Тогда верны неравенства
$$
\|f\|_{p,A}\ge \Bigl(\frac{\mu(A)}{c}\Bigr)^d (dp+1)^{-1/p} \|f\|_p, \text{ при $0<pd<1$},
$$
$$
\|f\|_{p,A}\ge \Bigl(\frac{\mu(A)}{c p d}\Bigr)^d \bigl(p+1/d)^{-1/p} \|f\|_p, \text{ при $pd\ge1$}.
$$

Отметим, что аналогичный результат верен для всякой системы функций, чьи сужения на отрезки удовлетворяют классическому
неравенству типа Ремеза. При этом достаточно заменить число $c$ на $3R^2$,
где $R$ --- это константа в классическом одномерном неравенстве Ремеза. Если для многочленов $A=4$, то для
тригонометрических многочленов подойдет $A=316$.

Теорема 1 является усилением следующего результата, см. \cite{AK}.

{\bf Теорема.} Пусть $f$ --- многочлен степени $d$ на локально выпуклом пространстве $E$,
причем $f$ является непрерывным многочленов или пределом по мере непрерывных многочленов.
Тогда найдется универсальная постоянная $C$, такая что для всякой логарифмически вогнутой меры $\mu$
имеет место неравенство
$$
\int_{A} |f|\ d\mu \ge \frac{\mu(A)^{d+1}}{(Cd)^{2d}} \|f\|_1
$$
для всякого измеримого множества $A$.

В цитированном результате величина постоянной содержит множитель $d^{2d}$, который
растет существенно медленнее полученного нами при $p=1$ множителя $d^{d+1}$.
Кроме того, в этой работе многочлен является числовым и на него наложены некоторые ограничения типа непрерывности,
в то время как в нашей работе накладывается лишь естественное ограничение измеримости.
В работе \cite{Ar} получено интегральное неравенство Ремеза
для произвольных измеримых полиномиальных отображений,
однако в этой работе константа имеет неявный вид (не ясно, как она зависит от
множества $A$ и степени $d$).

\section{Доказательство основного результата.}

Приведем неравенство Карбери---Райта.

{\bf Теорема.} Пусть $f$ --- это полиномиальное отображение степени $d$ в нормированное пространство $F$.
Тогда найдется универсальная постоянная $c$, такая что для всякой логарифмически вогнутой меры $\mu$
при $pd\ge1$ верно неравенство
$$
\mu(|f|\le t) \|f\|_p^{1/d} \le c t^{1/d} pd,
$$
а при $0<p<\frac{1}{d}$ верно неравенство
$$
\mu(|f|\le t) \|f\|_p^{1/d} \le c t^{1/d}.
$$

Аналогичное неравенство было получено Карбери и Райтом в статье \cite{CW} для конечномерных полиномов. Примерно в то же время
этот результат был получен Назаровым, Вольбергом и Содиным несколько другим методом в работе \cite{NVS}.
Одним из преимуществ этого метода явилась возможность перенести его на бесконечномерный случай (см. \cite{AK}).
В цитированной работе данное неравенство доказано в бесконечномерном случае для числовых многочленов.
Как перенести данное неравенство на случай полиномиальных отображений в банаховы пространства мы обсудим в следующем разделе.

\ \

{\bf Доказательство теоремы 1.} Без ограничения общности считаем, что $\|f\|_p=1$.
Обозначим $S=\bigl(\frac{\mu(A)}{cpd}\bigr)^d$ при $pd\ge1$ и $S=\bigl(\frac{\mu(A)}{c}\bigr)^d$ иначе. Тогда
$$
\int_A |f|^p \ d\mu = \int_{0}^{\infty} \mu(\{|f|^p\ge t\} \cap A) \ dt =
$$
$$
\int_{0}^{\infty} \mu(\{|f|\ge x\} \cap A)px^{p-1} \ dx\ge
\int_{0}^{\infty} \bigl(\mu(A) - \mu(\{|f|\le x\} \cap A)\bigr)px^{p-1} \ dx \ge
$$
$$
\int_{0}^{S} \bigl(\mu(A) - \mu(\{|f|\le x\} \cap A)\bigr)px^{p-1} \ dx \ge
\int_{0}^{S} \bigl(\mu(A) - \mu(\{|f|\le x\})\bigr)px^{p-1} \ dx.
$$
При $0<p<\frac{1}{d}$ последний интеграл оценивается как
$$
\int_{0}^{S} \bigl(\mu(A) - c x^{1/d}\bigr)px^{p-1} \ dx =
\mu(A)\Bigl(\frac{\mu(A)}{cpd}\Bigr)^{dp} - \frac{c p^2 d}{p+1/d} \Bigl(\frac{\mu(A)}{cpd}\Bigr)^{dp+1} =
$$
$$
\frac{\mu(A)^{dp+1}}{{c}^{dp}}\frac{1}{dp+1}.
$$
При $p>\frac{1}{d}$ оценка примет вид
$$
\int_{0}^{S} \bigl(\mu(A) - c pd x^{1/d}\bigr)px^{p-1} \ dx =
\mu(A)\Bigl(\frac{\mu(A)}{c}\Bigr)^{dp} - \frac{c p}{p+1/d} \Bigl(\frac{\mu(A)}{c}\Bigr)^{dp+1} =
$$
$$
\frac{\mu(A)^{dp+1}}{{cpd}^{dp}}\frac{d}{dp+1}.
$$

{\bf Замечание.} С точностью до некоторой константы в степени $d$ неравенство точное.
Рассмотрим для $p=1$ меру $\mu(dt)=e^{-t}I_{[0,\infty]}dt$.
$$
\|t^d\|=\int_0^{\infty} t^d e^{-t} = d! \ge \Bigl(\frac{d}{e}\Bigl)^d,
$$
$$
\int_{0}^{\varepsilon} t^d e^{-t} \le \frac{\varepsilon^{d+1}}{d+1}.
$$

{\bf Замечание.} При $p<0$ оценка величины $\|f\|_{L_p}$ через величину $\|f\|_{p,A}$ является тривиальной
даже для произвольной функции.
Действительно:
$$
\int_A \|f\|_F^p \le \int \|f\|_F^p,
$$
поэтому $\|f\|_p\le \bigl(\mu(A)\bigr)^{1/p} \|f\|_{p,A}$.
В то же время, в духе доказательства теоремы 1 можно получить неравенство, обратное тривиальному,
т.е. оценку сверху нормы сужения через норму по всему пространству,
однако в этом случае зависимость от меры множества и степени будет иметь несколько более сложный вид.

\section{Обобщение на случай полиномиальных отображений.}

Цель данного раздела заключается в том, чтобы доказать неравенство Карбери--Райта для норм полиномиальных отображений.
Сначала получим обобщение классического неравенства типа Ремеза на случай полиномиальных отображений.
Пусть $F$ --- нормированное пространство с нормой $\|\cdot\|_F$.
Тогда верно следующее неравенство.

{\bf Теорема 3.1.} Пусть $f$ --- полиномиальное отображение из $\mathbb{R}$ в $F$ степени $d$.
Пусть $\Delta$ --- некоторый отрезок на прямой, а $\lambda_{\Delta}$ --- это нормированная мера Лебега на отрезке $\Delta$.
Тогда для всякого измеримого множества $A\subset \Delta$ верно неравенство
$$
\sup_{\Delta} \|f\|_F \le \Bigl(\frac{4}{\lambda_{\Delta} (A)}\Bigr)^d \sup_{A} \|f\|_{F}.
$$
{\bf Доказательство.} Пусть $\ell_1, \ldots, \ell_n$ --- это непрерывные линейные функционалы на пространстве $F$.
Обозначим через $p_n$ полунорму на $F$, которая задается равенством $p_n(x)=\max\limits_{i=1,\ldots,n} |\ell_i(x)|$,
а также для четных чисел $k$ определим полунорму $p_{n,k}$ равенством $p_{n,k}(x) = (\frac{1}{n}\sum_{i=1}^n \ell_i (x)^k)^{1/k}$.
Заметим, что функция $(p_{n,k} (f))^k$ является многочленом степени $dk$, поэтому можно применить к нему
оригинальное неравенство Ремеза:
$$
\sup_{\Delta} (p_{n,k} (f))^k \le \Bigl(\frac{4}{\lambda_{\Delta} (A)}\Bigr)^{kd} \sup_{A} (p_{n,k} (f))^k.
$$
Извлекая корень степени $k$ получаем
$$
\sup_{\Delta} p_{n,k} (f) \le \Bigl(\frac{4}{\lambda_{\Delta} (A)}\Bigr)^{d} \sup_{A} p_{n,k} (f).
$$
Теперь заметим, что последовательность функций $p_{n,k} (f)$ монотонно растет по $k$, поэтому
в пределе можно получить
$$
\sup_{\Delta} p_{n} (f) \le \Bigl(\frac{4}{\lambda_{\Delta} (A)}\Bigr)^{d} \sup_{A} p_{n} (f).
$$
Теперь достаточно подобрать функционалы $\ell_1,\ldots,\ell_n,\ldots$ так,
чтобы левая и правая часть полученного неравенства при $n$ стремящемся к бесконечности
стремилась к, соответственно, левой и правой части желаемого неравенства.
Для этого заметим, что образ отрезка под действием многочлена $f(\Delta)$
является компактом в $F$. Поэтому достаточно показать,
что если $K$ --- это компакт в $F$, то можно найти
последовательность функционалов $\ell_i, i\in \mathbb{N}$, такую, что
$\|x\|_F=\sup_i |\ell_i(x)|$ для всякого $x\in K$. В качестве такой последовательности
подойдет набор опорных функционалов для точек $x_i$, образующих всюду плотное множество на данном компакте.

\ \

Остается понять, как вывести само неравенство Карбери--Райта.
Назаров, Вольберг и Содин в работе \cite{NVS} показали, как можно из классического неравенства Ремеза
для класса функций вывести неравенство Карбери--Райта. В бесконечномерном случае это
проделано в работе \cite{AK}. Таким образом после вывода классического неравенства Ремеза
не составляет труда вывести нужное неравенство Карбери--Райта.

%\section{Сравнение с предыдущими результатами}

\section{Дополнение}

Приведем одно утверждение, показывающее, что логарифмически вогнутые меры являются пределами
равномерных распределений на выпуклых множествах в слабой топологии.
Равномерным распределением на выпуклом множестве мы
называем нормированное сужение меры Лебега на данное множество.
При этом предполагается, что данное выпуклое тело ограничено, а также то,
что его аффинная оболочка является конечномерной гиперплоскостью.

Следующая теорема является аналогом задачи 2.7.50 из книги \cite{DM},
в которой дано подробное указание к решению.
В качестве пространства $X$ в этой задаче рассматриваются пространства $\ell^2$ и $\mathbb{R}^{\infty}$,
однако доказательство тривиальным образом переносится и на более общий случай.

{\bf Теорема 5.1.} Пусть $X$ --- локально выпуклое пространство, снабженное борелевской
логарифмически вогнутой мерой $\mu$.
Пусть также известно, что $X$ обладает базисом Шаудера.
Тогда найдется последовательность выпуклых множеств, равномерные
распределения на которых слабо сходятся к мере $\mu$.

%{\bf Доказательство.} Пусть $P_n$ --- это проекторы на первые $n$ базисных векторов.
%Заметим, что меры $\mu \circ P_n^{-1}$ слабо сходятся к мере $\mu$,
%где $\mu \circ P_n^{-1} (\cdot) = \mu(P_n^{-1}(\cdot))$. Действительно,
%для всякой ограниченной непрерывной функции $f$, поскольку $P_n(x)$ сходится к $x$, при $n$
%стремящемся к бесконечности имеем:
%$$
%\lim_{n\to \infty} \int f d\mu \circ P_n^{-1} = \lim_{n\to \infty} \int f(P_n) d\mu = \int f d\mu.
%$$
%В свою очередь, меры $\mu \circ P_n^{-1}$ --- это логарифмически вогнутые меры,
%сконцентрированные на конечномерных подпространствах.
%Широко известен факт о том, что всякая логарифмически вогнутая мера
%на конечномерном пространстве $L$ приближается
%проекциями на первую компоненту равномерных распределений на выпуклых телах $V$
%в пространствах $L\oplus \mathbb{R}^n$ (мера с плотностью $e^{-V(x)}, x\in L$ может быть приближена
%равномерными распределениями на телах $|s|_{\mathbb{R}^n} \le 1-\frac{V(x)}{n}, s\in \mathbb{R}^n$). Проекция же тела $V$ является %пределом равномерных распределений
%на телах $V_k$, которые получаются из него с помощью сужения по второй компоненте, т.е.
%$V_k = L_k(V)$, где $L_k\colon L\oplus \mathbb{R}^n, \ L_k(l+s)=l+\frac{s}{k}$. Теперь утверждение о наличии нужной последовательности
%следует из метризуемости слабой сходимости мер.

{\bf Замечание 5.2.} Легко видеть, что в предыдущей теореме можно заменить условие наличия базиса Шаудера
наличием стохастического базиса или просто наличием последовательности непрерывных операторов $P_n$
с конечномерными образами, таких что $P_n(x)$ сходится к $x$ $\mu$-п.в.

\end{document}